\documentclass[10pt]{amsart}
\usepackage{amsmath}
\usepackage{amscd}
\usepackage{amssymb}
\usepackage{amsthm}
\RequirePackage{filecontents}
\usepackage{fullpage}
\usepackage[numbers]{natbib}  
\usepackage[draft]{hyperref}

\newenvironment{psmallmatrix}
  {\left(\begin{smallmatrix}}
  {\end{smallmatrix}\right)}

\theoremstyle{plain}
\newtheorem{prop}{Proposition}
\newtheorem{thrm}[prop]{Theorem}

\newtheorem{cor}[prop]{Corollary}

\theoremstyle{definition}

\newtheorem{rem}[prop]{Remark}
\newtheorem{ex}[prop]{Example}

\title{Rankin-Cohen brackets and Serre derivatives as Poincar\'e series}
\author{Brandon Williams }

\address{Department of Mathematics \\ University of California \\ Berkeley, CA 94720}

\email{btw@math.berkeley.edu}

\begin{document}

\nocite{*}

\maketitle

\begin{abstract} We give expressions for the Serre derivatives of Eisenstein and Poincar\'e series as well as their Rankin-Cohen brackets with arbitrary modular forms in terms of the Poincar\'e averaging construction, and derive several identities for the Ramanujan tau function as applications.
\end{abstract}

\section{Introduction}

Let $k \in 2 \mathbb{Z}$, $k \ge 4.$ To any $q$-series $\phi(q) = \phi(e^{2\pi i \tau}) = \sum_{n=0}^{\infty} a_n q^n$ on the upper half-plane $\tau \in \mathbb{H}$ whose coefficients grow slowly enough, one can construct a \textbf{Poincar\'e series} $$\mathbb{P}_k(\phi;\tau) = \sum_{M \in \Gamma_{\infty} \backslash \Gamma} \phi |_k M (\tau) = \frac{1}{2} \sum_{c,d} \sum_{n=0}^{\infty} a_n (c \tau + d)^{-k} e^{2\pi i n \frac{a \tau + b}{c \tau + d}}$$ that converges absolutely and uniformly on compact subsets and defines a modular form of weight $k$. Here, the first sum is taken over cosets of $\Gamma = SL_2(\mathbb{Z})$ by the subgroup $\Gamma_{\infty}$ generated by $\pm \begin{pmatrix} 1 & 1 \\ 0 & 1 \end{pmatrix}$ and the second over all coprime integers $c,d \in \mathbb{Z}.$ As usual, $|_k$ is the Petersson slash operator of weight $k$. (More generally one can also construct Poincar\'e series that are not holomorphic in this way; see \cite{O}, section 8.3 for some applications.) \\

It is easy to show that every modular form $f$ (of weight $k \ge 4$) can be written as a Poincar\'e series $\mathbb{P}_k(\phi)$: because $f$ can always be written as a linear combination of the Eisenstein series $E_k = \mathbb{P}_k(1)$ and the Poincar\'e series of exponential type $P_{k,N} = \mathbb{P}_k(q^N)$ of various indices $N$. However, expressions found by this argument tend to be messy because the coefficients of $P_{k,N}$ are complicated series over Kloosterman sums and special values of Bessel functions (\cite{I}, section 3.2). The most reliable way to produce Poincar\'e series with manageable Fourier coefficients seems to be to start with seed functions $\phi(\tau)$ that already behave in a manageable way under the action of $SL_2(\mathbb{Z})$.

\begin{ex} When $\phi = 1$ (a modular form of weight $0$), we obtain the normalized Eisenstein series as mentioned above: $$\mathbb{P}_k(1;\tau) = E_k(\tau) = 1 - \frac{2k}{B_k} \sum_{n=1}^{\infty} \sigma_{k-1}(n) q^n, \; \; q = e^{2\pi i \tau}, \; \sigma_{k-1}(n) = \sum_{d|n} d^{k-1},$$ where $B_k$ is the $k$-th Bernoulli number. More generally, if $\phi$ is a modular form of any weight $k$ then expanding formally yields \begin{align*} \mathbb{P}_{k+l}(\phi;\tau) &= \sum_{M \in \Gamma_{\infty} \backslash \Gamma} (c \tau + d)^{-k-l} \phi \left( \frac{a \tau + b}{c \tau + d} \right) \\ &= \sum_{M \in \Gamma_{\infty} \backslash \Gamma} (c \tau + d)^{-k-l} (c \tau + d)^k \phi(\tau) \\ &= \phi(\tau) E_l(\tau), \end{align*} where $M$ is the coset of $ \begin{psmallmatrix} a & b \\ c & d \end{psmallmatrix}$; although the expression $\mathbb{P}_{k+l}(\phi)$ makes sense only when $l$ is sufficiently large compared to the growth of the coefficients of $\phi$. In recent work \cite{W} the author has considered the Poincar\'e series $\mathbb{P}_k(\vartheta)$ constructed from what are essentially weight $1/2$ theta functions $\vartheta$, which seem to be useful for computing with vector-valued modular forms for Weil representations; the details are somewhat more involved but this is related to the example above.
\end{ex}

The motivation of this note was to consider the Poincar\'e series $\mathbb{P}_k(\phi)$ when $\phi$ is a \textbf{quasimodular form}, a more general class of functions which includes modular forms, their derivatives of all orders, and the series $$E_2(\tau) = 1 - 24 \sum_{n=1}^{\infty} \sigma_1(n) q^n$$ (cf. \cite{Z1}, section 5.3). We find that one obtains Rankin-Cohen brackets and Serre derivatives (see section 2 below for their definitions) of Eisenstein series and Poincar\'e series essentially from such forms $\phi$:

\begin{thrm} For any modular form $f \in M_k$ and $l \in 2 \mathbb{N}$, $l \ge 4$, and $m,N \in \mathbb{N}_0$, with  $l \ge k+2$ if $f$ is not a cusp form, set $$\phi(\tau) = q^N \sum_{r=0}^m (-1)^r \binom{k+m-1}{m-r} \binom{l+m-1}{r} N^{m-r} D^r f(\tau);$$ then $$[f,P_{l,N}]_m= \mathbb{P}_{k+l+2m}(\phi).$$
\end{thrm}

Here $D = \frac{1}{2\pi i} \frac{d}{d \tau} = q \frac{d}{dq}$. Since $\mathbb{P}_{k+l+2m}(\phi)$ is modular by construction, and since $[-,-]_m$ is bilinear and $P_{l,N}$, $N \in \mathbb{N}_0$ span all modular forms, this gives another proof of the modularity of Rankin-Cohen brackets (at least for large $l$). This seed function $\phi$ is formally the Rankin-Cohen bracket $[f,q^N]_m$ where $q^N$ is treated like a modular form of weight $l$, so by linearity we see that Rankin-Cohen brackets and Poincar\'e averaging ``commute'' in the following sense:

\begin{cor} Let $f$ be a modular form of weight $k$ and let $\phi$ be a $q$-series whose coefficients grow sufficiently slowly that $\mathbb{P}_l(\phi;\tau)$ is well-defined, and denote by $[f,\phi]_m$ the formal result of the $m$-th Rankin-Cohen bracket where $\phi$ is treated like a modular form of weight $l$ (where $l \ge k+2$ if $f$ is not a cusp form). Then $$[f, \mathbb{P}_l(\phi)]_m = \mathbb{P}_{k+l+2m}([f,\phi]_m).$$
\end{cor}

 This expression simplifies considerably for the Eisenstein series: $$[f,E_l]_m = \mathbb{P}_{k+l+2m}(\phi)$$ for the function $$\phi = (-1)^m\binom{l+m-1}{m} D^m f.$$ An equivalent result in this case has appeared in section 5 of \cite{Z3} (in particular see Proposition 6). There may be particular interest in the case that $f$ itself is an Eisenstein or Poincar\'e series as expressions of a different nature for the Rankin-Cohen brackets of two Poincar\'e series are known (e.g. \cite{DS}, section 6). \\

\begin{thrm} For any $m,N \in \mathbb{N}_0$ and $l \in 2 \mathbb{N}$ with $l \ge 2m+2$, set $$\phi(\tau) = q^N \sum_{r=0}^m \binom{m}{r} \frac{(l+m-1)!}{(l+m-r-1)!} (-E_2(\tau) / 12)^r N^{m-r};$$ then the $m$-th order Serre derivative (in the sense of section 2) of $P_{l,N}$ is $$\vartheta^{[m]} P_{l,N} = \mathbb{P}_{l+2m}(\phi).$$
\end{thrm}

Similarly, this seed function $\phi$ is formally the $m$-th Serre derivative of $q^N$ if one pretends that $q^N$ is a modular form of weight $l$; by linearity we find that Serre derivatives also commute with Poincar\'e averaging:

\begin{cor} Let $\phi$ be a $q$-series whose coefficients grow sufficiently slowly that $\mathbb{P}_l(\phi;\tau)$ is well-defined, and denote by $\vartheta^{[m]} \phi$ the formal result of the $m$-th order Serre derivative where $\phi$ is treated like a modular form of weight $l$ (where $l \ge 2m+2$). Then $$\vartheta^{[m]} \mathbb{P}_l(\phi) = \mathbb{P}_{l+2m}(\vartheta^{[m]} \phi).$$
\end{cor}

As before, this simplifies for the Eisenstein series: $$\vartheta^{[m]} E_l = \mathbb{P}_{l+2m}(\phi)$$ for the function $$\phi= \frac{(l+m-1)!}{(-12)^m (l-1)!} E_2^m.$$ It is interesting to compare this to Theorem 2 which suggests that the Serre derivative (at least of the Eisenstein series) is analogous to a Rankin-Cohen bracket with $E_2$. Similar observations have been made before (e.g. \cite{G}, section 2). \\

By computing Rankin-Cohen brackets and Serre derivatives of $P_{l,N} = 0$ in weights $l \le 10$ we can obtain new proofs of Kumar's identity (\cite{K}, eq. (14)) $$\tau(m) = -\frac{20 m^{11}}{m - 5/6} \sum_{n=1}^{\infty} \frac{\sigma_1(n) \tau(m+n)}{(m+n)^{11}}$$ and Herrero's identity (\cite{H}, eq. (1)) $$\tau(m) = -240m^{11} \sum_{n=1}^{\infty} \frac{\sigma_3(n) \tau(m+n)}{(m+n)^{11}}$$ that express the Ramanujan tau function in terms of special values of a shifted $L$-series introduced by Kohnen \cite{Ko}, as well as four additional identities of this form. Namely we find \begin{align*} \tau(m) &= -\frac{14m^8}{m - 7/12} \sum_{n=1}^{\infty} \frac{\sigma_1(n) \tau(m+n)}{(m+n)^8} \\ &= -\frac{16m^9}{m - 2/3} \sum_{n=1}^{\infty} \frac{\sigma_1(n) \tau(m+n)}{(m+n)^9} \\ &= -\frac{18m^{10}}{m - 3/4} \sum_{n=1}^{\infty} \frac{\sigma_1(n) \tau(m+n)}{(m+n)^{10}} \\ &= -240m^{10} \sum_{n=1}^{\infty} \frac{\sigma_3(n)\tau(m+n)}{(m+n)^{10}}. \end{align*} Here $\tau(m)$ is Ramanujan's tau function, i.e. the coefficient of $q^m$ in $\Delta(\tau) = q \prod_{n=1}^{\infty} (1 - q^n)^{24}$.  We can also compute the values of these series with $m=0$. Based on numerical computations it seems reasonable to guess that there are no other identities of this type. The details are worked out in section 5.

\section{Background and notation}

Let $\mathbb{H} = \{\tau = x + iy: \, y > 0\}$ be the upper half-plane and let $\Gamma$ be the group $\Gamma = SL_2(\mathbb{Z})$, which acts on $\mathbb{H}$ by $\begin{psmallmatrix} a & b \\ c & d \end{psmallmatrix} \cdot \tau = \frac{a \tau + b}{c \tau + d}$. A \textbf{modular form of weight $k$} is a holomorphic function $f : \mathbb{H} \rightarrow \mathbb{C}$ which transforms under $\Gamma$ by $$f\left( \frac{a \tau + b}{c \tau + d} \right) = (c \tau + d)^k f(\tau), \; \; M = \begin{pmatrix} a & b \\ c & d \end{pmatrix} \in \Gamma, \; \tau \in \mathbb{H}$$ and whose Fourier expansion involves only non-negative exponents: $f(\tau) = \sum_{n=0}^{\infty} a_n q^n$, $q = e^{2\pi i \tau}$. We denote by $M_k$ the space of modular forms of weight $k$ and by $S_k$ the subspace of cusp forms (which in this context means $a_0 = 0$). \\

The \textbf{Rankin-Cohen brackets} are bilinear maps $$[\cdot,\cdot]_n : M_k \times M_l \rightarrow M_{k+l+2n},$$ \begin{equation}[f,g]_n = \sum_{j=0}^n (-1)^j \binom{k+n-1}{n-j} \binom{l+n-1}{j} D^j f D^{n-j} g,\end{equation} where $D^j f(\tau) = \frac{1}{(2\pi i)^j} \frac{d^j}{d \tau^j} f(\tau) = \frac{1}{(2\pi i)^j} f^{(j)}(\tau)$. If $f(\tau) = \sum_{n=0}^{\infty} a_n q^n$ is the Fourier expansion of $f$ then $$D^j f(\tau) = \sum_{n=0}^{\infty} a_n n^j q^n;$$ and in particular, the Rankin-Cohen brackets preserve integrality of Fourier coefficients. For example, the first few brackets are $$[f,g]_0 = fg, \; \; [f,g]_1 = k f \cdot Dg - lg \cdot Df,$$ $$[f,g]_2 =  \frac{k(k+1)}{2} f \cdot D^2 g - (k+1)(l+1) Df \cdot Dg + \frac{l(l+1)}{2} D^2 f \cdot g.$$ These can be characterized as the unique (up to scale) bilinear differential operators of degree $2n$ that preserve modularity (see for example the second proof in section 1 of \cite{Z2}). \\

The Serre derivatives $\vartheta^{[n]}$ following \cite{Z1}, section 5.1 are maps $M_k \rightarrow M_{k+2n}$ defined recursively by $$\vartheta^{[0]} f = f, \; \; \vartheta^{[1]} f = \vartheta f = Df- \frac{k}{12} E_2 f,$$ and $$\vartheta^{[n+1]} f = \vartheta \vartheta^{[n]} f - \frac{n(k+n-1)}{144} E_4f, \; \; n \ge 1.$$ (In particular $\vartheta^{[n]}$ is \emph{not} simply the $n$-th iterate of $\vartheta$.) These functions are given in closed form by \begin{equation} \vartheta^{[n]} f(\tau) = \sum_{r=0}^n \binom{n}{r} \frac{(k+n-1)!}{(k+r-1)!} (-E_2(\tau) / 12)^{n-r} D^r f(\tau),\end{equation} as one can prove by induction or by inverting equation 65 of \cite{Z1} (section 5.2).

\section{Poincar\'e series}

\begin{rem} A sufficient criterion for the series $$\mathbb{P}_k(\phi;\tau) = \sum_{c,d} (c \tau + d)^{-k} \phi\Big( \frac{a \tau + b}{c \tau + d} \Big), \; \; \phi(\tau) = \sum_{n=0}^{\infty} a_n q^n$$ to converge absolutely and locally uniformly is for the coefficients of $\phi$ to satisfy the bound $a_n = O(n^l)$ where $l = \frac{k}{2} - 2 - \varepsilon$ for some $\varepsilon > 0$. To see this, note that $\binom{n+l}{l}$ is also $O(n^l)$, so we can bound $$\Big| \phi\Big( \frac{a \tau + b}{c \tau + d} \Big) \Big| \ll \sum_{n=0}^{\infty} \binom{n+l}{l} e^{-2\pi n \frac{1}{|c \tau + d|^2}} = \Big( 1 - e^{-\frac{2\pi}{|c \tau + d|^2}} \Big)^{-l-1}$$ up to a constant multiple. Since $(1 - e^{-x})^{-1} < x^{-1-\delta}$ for any fixed (small enough) $\delta > 0$ and all small enough $x > 0$, we can then bound $$\sum_{c,d} \Big| (c \tau + d)^{-k} \phi\Big( \frac{a \tau + b}{c \tau + d} \Big) \Big| \ll \sum_{c,d} |c \tau + d|^{-k + 2(l+1)(1 + \delta)} < \sum_{c,d} |c \tau + d|^{-2}.$$
\end{rem}

\begin{rem} Given a $q$-series $\phi(\tau) = \sum_{n=0}^{\infty} a_n q^n$, one can also consider the series $$\mathbb{P}_k'(\phi) = a_0 E_k + \sum_{n=1}^{\infty} a_n P_{k,n}$$ which generally has better convergence properties than the sum $\mathbb{P}_k(\phi)$ over cosets $\Gamma_{\infty} \backslash \Gamma$. Since $S_k$ is finite-dimensional, the convergence of $\sum_{n=1}^{\infty} a_n P_{k,n}$ to a cusp form in any sense is equivalent to the convergence of the series $$\sum_{n=1}^{\infty} a_n \langle f, P_{k,n} \rangle = \frac{(k-2)!}{(4\pi)^{k-1}} \sum_{n=1}^{\infty} \frac{a_n b_n}{n^{k-1}}$$ for every cusp form $f(\tau) = \sum_{n=1}^{\infty} b_n q^n \in S_k$. The Deligne bound $b_n = O(n^{(k-1)/2 + \varepsilon})$ implies that this is satisfied when the slightly weaker bound $a_n = O(n^{k/2 - 3/2 - \varepsilon})$ holds. It is clear that $\mathbb{P}_k'(\phi) = \mathbb{P}_k(\phi)$ whenever the latter series converges, so we will refer to both of these series by $\mathbb{P}_k(\phi;\tau)$ in what follows.
\end{rem}

\section{Proofs}

\begin{proof}[Proof of Theorem 2]

The coefficients $a_n$ of any modular form of weight $k$ satisfy the bound $a_n = O(n^{k-1 + \varepsilon})$ for any $\varepsilon > 0$, while cusp forms satisfy the Deligne bound $a_n = O(n^{(k-1)/2 + \varepsilon})$. In particular, the coefficients of $$\phi(\tau) = q^N \sum_{r=0}^m (-1)^r \binom{k+m-1}{n-r} \binom{l+m-1}{r} N^{m-r} D^r f(\tau)$$ always satisfy the bound $O(n^{k+m-1 + \varepsilon})$, and our growth condition (of Remark 7), $$k+m-1+ \varepsilon \le \frac{k+l+2m}{2} - 3/2 - \varepsilon$$ becomes $k \le l - 1 - 2 \varepsilon$ and therefore (since $k,l \in 2 \mathbb{Z}$) $l \ge k+2$; while for cusp forms we instead require $$\frac{k-1}{2} + m + \varepsilon \le \frac{k+l+2m}{2} - 3/2 - \varepsilon,$$ or equivalently $2 \le l - 2 \varepsilon$ which is always satisfied for small enough $\varepsilon$. \\

Suppose first that the series $\mathbb{P}(\phi;\tau)$ over cosets $\Gamma_{\infty} \backslash \Gamma$ converges normally. Repeatedly differentiating the equation $$f\Big( \frac{a \tau + b}{c \tau + d} \Big) = (c \tau + d)^k f(\tau)$$ yields $$f^{(m)}\Big( \frac{a \tau + b}{c \tau + d} \Big) = \sum_{r=0}^m \binom{m}{r} \frac{(k+m-1)!}{(k+r-1)!} c^{m-r} (c \tau + d)^{k+m+r} f^{(r)}(\tau),$$ as one can prove by induction or derive directly by considering the action of $SL_2(\mathbb{Z})$ on $\tau$ in the generating series $$\sum_{m=0}^{\infty} f^{(m)}(\tau) \frac{w^m}{m!} = f(\tau + w),$$ for $|w|$ sufficiently small. By another induction argument one finds the similar formula \begin{equation} \frac{d^m}{d \tau^m} \Big( (c \tau + d)^{-k} e^{2\pi i N \tau}\Big) = \sum_{r=0}^m \binom{m}{r} \frac{(k+m-1)!}{(k+r-1)!} (-c)^{m-r} (2\pi i N)^r (c \tau + d)^{-k-m-r} e^{2\pi i N \tau} \end{equation} for any $N \in \mathbb{N}_0$. \\

Let $(a)_m = \frac{(a+m-1)!}{(a-1)!} = a \cdot (a+1) \cdot ... \cdot (a+m-1)$ denote the Pochhammer symbol. Then \begin{align*} &\quad \sum_{M \in \Gamma_{\infty} \backslash \Gamma} \Big[ q^N \sum_{r=0}^m (-1)^r \binom{k+m-1}{m-r} \binom{l+m-1}{r} N^{m-r} D^r f(\tau) \Big] \Big|_{k+l+2m} M (\tau) \\ &= \sum_{M \in \Gamma_{\infty} \backslash \Gamma} \sum_{r=0}^m \sum_{j=0}^r \Big[ \binom{k+m-1}{m-r} \binom{l+m-1}{r} \binom{r}{j} N^{m-r} \times \\ &\quad \quad = (-2\pi i)^{-r} c^{r-j} (k+j)_{r-j} (c \tau + d)^{j+r-l-2m} e^{2\pi i N \frac{a \tau + b}{c \tau + d}} f^{(j)}(\tau) \Big] \\ &=(-1)^m \sum_{M \in \Gamma_{\infty} \backslash \Gamma} \sum_{j=0}^m (-1)^j f^{(j)}(\tau) \sum_{r=0}^{m-j}\Big[ (2\pi i N)^r \binom{k+m-1}{r} \binom{l+m-1}{m-r} \binom{m-r}{j} (k+j)_{m-j-r} \times \\ &\quad\quad\quad\quad\quad\quad\quad\quad \times  (-c)^{m-r-j} (c \tau + d)^{j-l-m-r} e^{2\pi i N \frac{a \tau + b}{c \tau + d}} \Big],\end{align*} where we have replaced $r$ by $m-r$ in the second equality. Since \begin{align*} &\quad \binom{k+m-1}{r} \binom{l+m-r}{m-r} \binom{m-r}{j} (k+j)_{m-j-r} \\ &= \frac{(k+m-1)! (l+m-1)! (m-r)! (k+m-r-1)!}{r!(k+m-r-1)!(m-r)!(l+r-1)!j!(m-r-j)!(k+j-1)!} \\ &= \binom{k+m-1}{m-j} \binom{l+m-1}{j} \binom{m-j}{r} (l+r)_{m-j-r}, \end{align*} as we see by replacing $\frac{(m-r)! (k+m-r-1)!}{(m-r)!(k+m-r-1)!}$ by $\frac{(m-j)! (l+m-j-1)!}{(m-j)!(l+m-j-1)!}$ in the above expression, this equals \begin{align*} &\quad (2\pi i)^{-m} \sum_{j=0}^m \Big[ (-1)^j f^{(j)}(\tau) \binom{k+m-1}{m-j} \binom{l+m-1}{j} \times \\ &\quad\quad\quad\quad \times \sum_{M \in \Gamma_{\infty} \backslash \Gamma} \sum_{r=0}^{m-j} (2\pi i N)^r \binom{m-j}{r} (l+r)_{m-j-r} (-c)^{m-r-j} (c \tau + d)^{j-l-m-r} e^{2\pi i N \frac{a\tau +b}{c \tau + d}}\Big] \\ &= \sum_{j=0}^m (-1)^j D^j f(\tau) \binom{k+m-1}{m-j} \binom{l+m-1}{j} \sum_{M \in \Gamma_{\infty} \backslash \Gamma} D^{m-j} \Big( (c \tau + d)^{-l} e^{2\pi i N \frac{a \tau + b}{c \tau + d}} \Big) \\ &= \sum_{j=0}^m (-1)^j \binom{k+m-1}{m-j} \binom{l+m-1}{j} D^j f(\tau) D^{m-j} P_{l,N}(\tau) \\ &= [f,P_{l,N}]_m(\tau),\end{align*} the last equality by definition (equation (1)), and the third-to-last equality using equation (3). \\

When $\phi$ satisfies the weaker growth condition, we can include a convergence factor $(c \overline{\tau} + d)^{-s}$ into the argument above (which is ignored by the operator $D$) to see that, if $\phi(\tau) = a_0 + a_1 q + a_2 q^2 + ...$, then
\begin{align*} &\quad a_0 E_k(\tau;s) + a_1 P_{k,1}(\tau;s) + a_2 P_{k,2}(\tau;s) + ...  \\  &= \sum_{M \in \Gamma_{\infty} \backslash \Gamma} (c \overline{\tau} + d)^{-s} \times \phi(\tau) \Big|_{k+l+2m} M \\ &= \sum_{j=0}^m (-1)^j \binom{k+m-1}{m-j} \binom{l+m-1}{j} D^j f(\tau) D^{m-j} \Big( \sum_{M \in \Gamma_{\infty} \backslash \Gamma} (c \overline{\tau} + d)^{-s} \times q^N \Big|_l M \Big)\end{align*} when $\mathrm{Re}[s]$ is sufficiently large, and $E_k(\tau;s)$ and $P_{k,N}(\tau;s)$ denote the deformed series $$E_k(\tau;s) = \frac{1}{2}\sum_{c,d} \frac{1}{(c \tau + d)^k (c \overline{\tau} + d)^s}, \; \; P_{k,N}(\tau;s) = \frac{1}{2} \sum_{c,d} \frac{e^{2\pi i N \frac{a \tau + b}{c \tau + d}}}{(c \tau + d)^k (c \overline{\tau} + d)^s}.$$ The claim follows by analytic continuation to $s=0$.

\end{proof}

\begin{proof}[Proof of Theorem 4] The condition $l \ge 2m+2$ makes the Fourier coefficients of $\phi$ grow sufficiently slowly: the $n$-th coefficient of $E_2^m$ is $O(n^{2m-1 + \varepsilon})$ for any $\varepsilon > 0$, so the growth condition $$2m-1+\varepsilon \le \frac{l+2m}{2} - 3/2 - \varepsilon$$ of Remark 7 is satisfied for all $l \ge 2m + 2$. \\

Using the transformation law $$E_2 \Big( \frac{a \tau + b}{c \tau + d} \Big) = (c \tau + d)^2 E_2(\tau) + \frac{6}{\pi i} c (c \tau + d),$$ it follows that $$E_2 \Big( \frac{a \tau + b}{c \tau + d} \Big)^m = \sum_{r=0}^m \binom{m}{r} (c \tau + d)^{m+r} c^{m-r} \Big( \frac{12}{2\pi i} \Big)^{m-r} E_2(\tau)^r$$ for all $m \in \mathbb{N}$. Therefore, with $$\phi(\tau) = q^N \sum_{r=0}^m \binom{m}{r} \frac{(l+m-1)!}{(l+m-r-1)!} (-E_2(\tau) / 12)^r N^{m-r},$$ ignoring convergence issues for now, we find

\begin{align*} &\quad \sum_{M \in \Gamma_{\infty} \backslash \Gamma} \phi \Big|_{l+2m} M (\tau) \\ &= \sum_{r=0}^m (-12)^{-r} N^{m-r} \binom{m}{r} \frac{(l+m-1)!}{(l+m-r-1)!} \sum_M \sum_{j=0}^r \binom{r}{j} c^{r-j} (c \tau + d)^{r+j}  (12 / 2\pi i)^{r-j} E_2(\tau)^j e^{2\pi i N \frac{a \tau + b}{c \tau + d}} \\ &= \sum_{j=0}^m \sum_{r=j}^m (-12)^{-r} N^{m-r} \binom{m}{r} \frac{(l+m-1)!}{(l+m-r-1)!} \binom{r}{j} (12 / 2\pi i)^{r-j} E_2(\tau)^j \sum_M \Big[ c^{r-j} (c \tau + d)^{r+j-l-2m} e^{2\pi i N \frac{a \tau + b}{c \tau + d}} \Big] \\ &= \sum_{j=0}^m E_2(\tau)^j \sum_{r=0}^{m-j} (-12)^{r-m} N^r \binom{m}{r} \frac{(l+m-1)!}{(l+r-1)!} \binom{m-r}{j} (12 / 2\pi i)^{m-r-j} \sum_M c^{m-r-j} (c \tau + d)^{j-l-m-r} e^{2\pi i N \frac{ a \tau + b}{c \tau + d}}, \end{align*} where in the last line we replaced $r$ by $m-r$. Since \begin{align*} &\quad (-12)^{r-m} N^r \binom{m}{r} \frac{(l+m-1)!}{(l+r-1)!} \binom{m-r}{j} (-12 / 2\pi i)^{m-r-j} \\ &= (2\pi i)^{-m} \binom{m}{j} \frac{(l+m-1)!}{(l+m-j-1)!} (-2\pi i / 12)^j \binom{m-j}{r} \frac{(l+m-j-1)!}{(l+r-1)!} (2\pi i N)^r, \end{align*} as one can see by expanding both sides of this equation, the expression above equals \begin{align*} &\quad (2\pi i)^{-m} \sum_{j=0}^m E_2(\tau)^j \sum_{r=0}^{m-j} \Big[ \binom{m}{j} \frac{(l+m-1)!}{(l+m-j-1)!} (-2\pi i / 12)^j \binom{m-j}{r} \frac{(l+m-j-1)!}{(l+r-1)!} \times  \\ &\quad\quad\quad\quad\quad\quad\quad\quad \times \sum_M (-c)^{m-j-r} (c \tau + d)^{j-l-m-r} (2\pi i N)^r e^{2\pi i N \frac{a \tau + b}{c \tau + d}} \Big] \\ &= \sum_{j=0}^m \binom{m}{j} \frac{(l+m-1)!}{(l+m-j-1)!} \Big( -E_2(\tau) / 12 \Big)^j D^{m-j} P_{l,N}(\tau) \\ &= \vartheta^{[m]} P_{l,N}(\tau), \end{align*} using equation (2) from section 2. Convergence issues can be resolved by including the factor $(c \overline{\tau} + d)^{-s}$ as in the proof of Theorem 2.
\end{proof}

\section{Examples involving Ramanujan's tau function}

In weight $12$, the space $S_k$ of cusp forms is one-dimensional and therefore all Poincar\'e series are multiples of the discriminant $\Delta(\tau) = \sum_{n=1}^{\infty} \tau(n) q^n$; we find this multiple by writing $P_{k,m} = \lambda_m \Delta$ and using $\lambda_m \langle \Delta, \Delta \rangle = \langle \Delta, P_{k,m} \rangle = \tau(m) \frac{10!}{(4\pi m)^{11}},$ such that $$P_{k,m} = \frac{10! \cdot \tau(m)}{(4\pi m)^{11} \langle \Delta, \Delta \rangle} \Delta.$$ We can form the Poincar\'e series $\mathbb{P}_{12}(\phi)$ from any $q$-series $\phi(\tau) = \sum_{n=0}^{\infty} a_n q^n$ with $a_n = O(n^{9/2 - \varepsilon})$. This includes the $q$-series $E_2$ and $E_4$ and some of their derivatives. Applying Theorems $2$ and $4$ together with the vanishing of cusp forms in weight $\le 10$ gives identities involving $\tau(n)$. (Similar arguments can be used to derive identities for the coefficients of the normalized cusp forms of weights $16,18,20,22,26$.)

\begin{ex} By Theorem 4, $$0 = \vartheta P_{10,m} = \mathbb{P}_{12}\Big[ q^m \Big( m - \frac{5}{6}E_2 \Big) \Big] = (m-5/6) P_{12,m} + (-5/6) \cdot (-24) \sum_{n=1}^{\infty} \sigma_1(n) P_{12,m+n},$$ so we recover Kumar's identity $$\tau(m) = -\frac{20 m^{11}}{m - 5/6} \sum_{n=1}^{\infty} \frac{\tau(m+n) \sigma_1(n)}{(m+n)^{11}}.$$
\end{ex}

\begin{ex} By Theorem 2, $$0 = P_{8,m}E_4 = \mathbb{P}_{12}(q^m E_4) = P_{12,m} - 240 \sum_{n=1}^{\infty} \sigma_3(n) P_{12,m+n},$$ which yields Herrero's identity $$\tau(m) = -240 m^{11} \sum_{n=1}^{\infty} \frac{\sigma_3(n) \tau(m+n)}{(m+n)^{11}}.$$
\end{ex}


\begin{ex} By Theorem 2, $$0 = [E_4,P_{6,m}]_1 = \mathbb{P}_{12} \Big( 4m E_4 + 6 DE_4\Big) = 4m P_{12,m} + 240 \sum_{n=1}^{\infty} (4m + 6n) \sigma_3(n) P_{12,m+n},$$ which implies $$\tau(m) = -60 m^{10} \sum_{n=1}^{\infty} \frac{(4m + 6n) \sigma_3(n) \tau(m+n)}{(m+n)^{11}}.$$ Together with the previous identity this implies $$\tau(m) = -240 m^{10} \sum_{n=1}^{\infty} \frac{(m+n) \sigma_3(n) \tau(m+n)}{(m+n)^{11}} = -240m^{10} \sum_{n=1}^{\infty} \frac{\sigma_3(n) \tau(m+n)}{(m+n)^{10}}.$$
\end{ex}

\begin{ex} By Theorem 4, $$0 = \vartheta^{[2]} P_{8,m} = \mathbb{P}_{12}\Big( q^m (m^2 - (3/2)m E_2(\tau) + (1/2) E_2(\tau)^2) \Big),$$ where by Ramanujan's equation $D E_2 = \frac{1}{12} (E_2^2 - E_4)$ the coefficient of $q^n$ in $E_2(\tau)^2$ is $240\sigma_3(n) - 288n \sigma_1(n)$. Therefore we find $$0 = \left( m^2 - \frac{3}{2}m + \frac{1}{2} \right) P_{12,m} + \sum_{n=1}^{\infty} \Big( 36m \sigma_1(n) + 120 \sigma_3(n) - 144 n \sigma_1(n) \Big) P_{12,m+n}$$ and therefore $$(2m^2 - 3m + 1)\tau(m) = -24m^{11} \sum_{n=1}^{\infty} \frac{((3m-12n) \sigma_1(n) + 10 \sigma_3(n)) \tau(m+n)}{(m+n)^{11}}, \;\; m \in \mathbb{N}.$$ Combining this with the previous examples, we find $$\tau(m) = -180m^9 \sum_{n=1}^{\infty} \frac{n \sigma_1(n) \tau(m+n)}{(m+n)^{11}}$$ and therefore $$\tau(m) = -\frac{18m^{10}}{m - 3/4} \sum_{n=1}^{\infty} \frac{\sigma_1(n) \tau(m+n)}{(m+n)^{10}}.$$ 
\end{ex}

\begin{ex} It is not valid to form the Poincar\'e series $\mathbb{P}_{12}(\phi)$ with either $\phi = E_2^3$ or $E_6$, because their Fourier coefficients grow too quickly; however, their difference $E_2^3 - E_6 = 9 DE_4 +72 D^2 E_2$ has coefficients that satisfy the required bound $O(n^{9/2-\varepsilon})$. We use \begin{align*} 0 &= \vartheta^{[3]} P_{6,m} + \frac{7}{36} P_{6,m} E_6 \\ &= \mathbb{P}_{12}\Big( q^m (m^3 - 2m^2 E_2 + (7/6) m E_2^2 - (7/36) (E_2^3 - E_6)) \Big) \\ &= \left(m^3 - 2m^2 + \frac{7}{6}m \right) P_{12,m} + \sum_{n=1}^{\infty} \Big[ (48m^2 - 336mn + 336n^2) \sigma_1(n) + (280m - 420n) \Big] P_{12,m+n} \end{align*} to obtain $$\tau(m) = -720m^8 \sum_{n=1}^{\infty} \frac{n^2 \sigma_1(n) \tau(m+n)}{(m+n)^{11}},$$ and combining this with the previous examples, $$\tau(m) = -\frac{16m^9}{m - 2/3} \sum_{n=1}^{\infty} \frac{\sigma_1(n) \tau(m+n)}{(m+n)^9}.$$ Similarly, by expressing $D^3 E_2$ in terms of powers of $E_2$ and derivatives of modular forms one obtains the formula $$\tau(m) = -\frac{14m^8}{m - 7/12} \sum_{n=1}^{\infty} \frac{\sigma_1(n) \tau(m+n)}{(m+n)^8}.$$
\end{ex}

\begin{rem} In particular, for any $m \in \mathbb{N}$ the values of the $L$-series $\sum_{n=1}^{\infty} \frac{\sigma_1(n) \tau(m+n)}{(m+n)^s}$ at $s=8,9,10,11$ and of $\sum_{n=1}^{\infty} \frac{\sigma_3(n) \tau(m+n)}{(m+n)^s}$ at $s=10,11$ are rational numbers, and Lehmer's conjecture that $\tau(n)$ is never zero is equivalent to the non-vanishing of any of these $L$-values. Computing these $L$-series at other integers $s$ numerically does not seem to yield rational numbers. In any case, the methods of this note do not apply to other values of $s$.
\end{rem}

We can also evaluate the values of these $L$-series with $m=0$ by a similar argument. Comparing $$\vartheta E_{10} = -\frac{5}{6} - 24q - ... = -\frac{5}{6} E_{12} + \frac{38016}{691} \Delta$$ with the result of Theorem 4, $$\vartheta E_{10} = -\frac{5}{6} E_{12} + 20 \sum_{n=1}^{\infty} \sigma_1(n) P_{12,n}$$ we find $$\tau(m) = \frac{20 \cdot 691}{38016} \sum_{n=1}^{\infty} \tau(n) \tau(m) \sigma_1(n) \cdot \frac{10!}{\langle \Delta, \Delta \rangle \cdot (4\pi n)^{11}},$$ i.e. $$\sum_{n=1}^{\infty} \frac{\tau(n) \sigma_1(n)}{n^{11}} = \frac{2^{19} \cdot 11}{3 \cdot 5^3 \cdot 7 \cdot 691} \pi^{11} \langle \Delta, \Delta \rangle \approx 0.968.$$ Here, the Petersson norm-square of $\Delta$ to $18$ decimal places is $$\langle \Delta, \Delta \rangle \approx 1.03536205680 \times 10^{-6}$$ which can be computed using PARI/GP. \\

Similarly, comparing $E_8 E_4 = 1 + 720q + ... = E_{12} + \frac{432000}{691} \Delta$ with $$E_8 E_4 = \mathbb{P}_{12}(E_4) = E_{12} + 240 \sum_{n=1}^{\infty} \sigma_3(n) P_{12,n},$$ we find $$\tau(m) = \frac{240 \cdot 691}{432000} \sum_{n=1}^{\infty} \tau(n) \tau(m) \sigma_3(n) \cdot \frac{10!}{\langle \Delta, \Delta \rangle \cdot (4\pi n)^{11}},$$ i.e. $$\sum_{n=1}^{\infty} \frac{\tau(n) \sigma_3(n)}{n^{11}} = \frac{2^{17}}{3^2 \cdot 7 \cdot 691} \pi^{11} \langle \Delta, \Delta \rangle \approx 0.917.$$

With similar arguments applied to \begin{align*} -3456 \Delta &= [E_4,E_6]_1 = -6 \mathbb{P}_{12}(DE_4), \\ \frac{1}{2} E_{12} - \frac{49344}{691} \Delta &= \vartheta^{[2]} E_8 = \frac{1}{2} \mathbb{P}_{12}(E_2^2), \\ -168 \Delta &= \vartheta^{[3]} E_6 + \frac{7}{36} E_6^2 = \frac{7}{36} \mathbb{P}( E_6 - E_2^3), \end{align*} and $$-600 \Delta = \vartheta^{[4]}E_4 - \frac{35}{864} E_4 E_8 - \frac{7}{40} [E_4,E_4]_2 + \frac{35}{432} [E_6,E_4]_1 = \frac{35}{3} \mathbb{P}_{12}(D^3 E_2),$$ one can compute the values \begin{align*} \sum_{n=1}^{\infty} \frac{\tau(n) \sigma_3(n)}{n^{10}} &= \frac{2^{16}}{3^3 \cdot 5^3 \cdot 7} \pi^{11} \langle \Delta, \Delta \rangle \approx 0.845, \\ \sum_{n=1}^{\infty} \frac{\tau(n) \sigma_1(n)}{n^{10}} &= \frac{2^{17}}{3^5 \cdot 5^2 \cdot 7} \pi^{11} \langle \Delta, \Delta \rangle  \approx 0.939, \\ \sum_{n=1}^{\infty} \frac{\tau(n) \sigma_1(n)}{n^9} &= \frac{2^{13}}{3^4 \cdot 5 \cdot 7} \pi^{11} \langle \Delta,\Delta \rangle \approx 0.880, \\ \sum_{n=1}^{\infty} \frac{\tau(n) \sigma_1(n)}{n^8} &= \frac{2^{14}}{3^3 \cdot 5 \cdot 7^2} \pi^{11} \langle \Delta, \Delta \rangle \approx 0.754. \end{align*} Unlike the $L$-values of examples $8$ through $12$, none of these are expected to be rational. \\

\textbf{Acknowledgments:} I thank the reviewers for pointing out a mistake in Remark 6 in an earlier version of this note, and also for several suggestions that improved the exposition.

\bibliographystyle{plainnat}
\bibliography{\jobname}

\begin{thebibliography}{11}
\providecommand{\natexlab}[1]{#1}
\providecommand{\url}[1]{\texttt{#1}}
\expandafter\ifx\csname urlstyle\endcsname\relax
  \providecommand{\doi}[1]{doi: #1}\else
  \providecommand{\doi}{doi: \begingroup \urlstyle{rm}\Url}\fi

\bibitem[Diamantis and O'Sullivan(2013)]{DS}
Nikolaos Diamantis and Cormac O'Sullivan.
\newblock Kernels for products of {$L$}-functions.
\newblock \emph{Algebra Number Theory}, 7\penalty0 (8):\penalty0 1883--1917,
  2013.
\newblock ISSN 1937-0652.
\newblock URL \url{https://doi.org/10.2140/ant.2013.7.1883}.

\bibitem[Gouv\^ea(1997)]{G}
Fernando Gouv\^ea.
\newblock Non-ordinary primes: a story.
\newblock \emph{Experiment. Math.}, 6\penalty0 (3):\penalty0 195--205, 1997.
\newblock ISSN 1058-6458.
\newblock URL \url{http://projecteuclid.org/euclid.em/1047920420}.

\bibitem[Herrero(2015)]{H}
Sebasti\'an~Daniel Herrero.
\newblock The adjoint of some linear maps constructed with the {R}ankin-{C}ohen
  brackets.
\newblock \emph{Ramanujan J.}, 36\penalty0 (3):\penalty0 529--536, 2015.
\newblock ISSN 1382-4090.
\newblock URL \url{https://doi.org/10.1007/s11139-013-9536-5}.

\bibitem[Iwaniec(1997)]{I}
Henryk Iwaniec.
\newblock \emph{Topics in classical automorphic forms}, volume~17 of
  \emph{Graduate Studies in Mathematics}.
\newblock American Mathematical Society, Providence, RI, 1997.
\newblock ISBN 0-8218-0777-3.
\newblock URL \url{https://doi.org/10.1090/gsm/017}.

\bibitem[Kohnen(1991)]{Ko}
Winfried Kohnen.
\newblock Cusp forms and special values of certain {D}irichlet series.
\newblock \emph{Math. Z.}, 207\penalty0 (4):\penalty0 657--660, 1991.
\newblock ISSN 0025-5874.
\newblock URL \url{https://doi.org/10.1007/BF02571414}.

\bibitem[Kumar(2017)]{K}
Arvind Kumar.
\newblock The adjoint map of the {S}erre derivative and special values of
  shifted {D}irichlet series.
\newblock \emph{J. Number Theory}, 177:\penalty0 516--527, 2017.
\newblock ISSN 0022-314X.
\newblock URL \url{https://doi.org/10.1016/j.jnt.2017.01.011}.

\bibitem[Ono(2009)]{O}
Ken Ono.
\newblock Unearthing the visions of a master: harmonic {M}aass forms and number
  theory.
\newblock In \emph{Current developments in mathematics, 2008}, pages 347--454.
  Int. Press, Somerville, MA, 2009.

\bibitem[Williams()]{W}
Brandon Williams.
\newblock Poincar\'e square series for the {W}eil representation.
\newblock \emph{Ramanujan J.}, (in press).
\newblock \doi{10.1007/s11139-017-9986-2}.
\newblock URL \url{https://doi.org/10.1007/s11139-017-9986-2}.

\bibitem[Zagier(1977)]{Z3}
Don Zagier.
\newblock Modular forms whose {F}ourier coefficients involve zeta-functions of
  quadratic fields.
\newblock pages 105--169. Lecture Notes in Math., Vol. 627, 1977.

\bibitem[Zagier(1994)]{Z2}
Don Zagier.
\newblock Modular forms and differential operators.
\newblock \emph{Proc. Indian Acad. Sci. Math. Sci.}, 104\penalty0 (1):\penalty0
  57--75, 1994.
\newblock ISSN 0253-4142.
\newblock URL \url{https://doi.org/10.1007/BF02830874}.
\newblock K. G. Ramanathan memorial issue.

\bibitem[Zagier(2008)]{Z1}
Don Zagier.
\newblock Elliptic modular forms and their applications.
\newblock In \emph{The 1-2-3 of modular forms}, Universitext, pages 1--103.
  Springer, Berlin, 2008.
\newblock URL \url{https://doi.org/10.1007/978-3-540-74119-0_1}.

\end{thebibliography}

\end{document}